\pdfoutput=1
\RequirePackage{ifpdf}
\ifpdf 
\documentclass[pdftex]{sigma}
\else
\documentclass{sigma}
\fi

\usepackage[all]{xy}

\DeclareMathOperator{\sgn}{sgn}

\begin{document}

\newcommand{\arXivNumber}{1406.1072}

\allowdisplaybreaks

\renewcommand{\thefootnote}{$\star$}

\renewcommand{\PaperNumber}{089}

\FirstPageHeading

\ShortArticleName{Maximal Green Sequences of Exceptional Finite Mutation Type Quivers}

\ArticleName{Maximal Green Sequences\\
of Exceptional Finite Mutation Type Quivers\footnote{This paper is a~contribution to the Special Issue on New Directions
in Lie Theory.
The full collection is available at
\href{http://www.emis.de/journals/SIGMA/LieTheory2014.html}{http://www.emis.de/journals/SIGMA/LieTheory2014.html}}}

\Author{Ahmet I.~SEVEN}

\AuthorNameForHeading{A.I.~Seven}

\Address{Middle East Technical University, Department of Mathematics, 06800, Ankara, Turkey}
\Email{\href{mailto:aseven@metu.edu.tr}{aseven@metu.edu.tr}}

\ArticleDates{Received June 18, 2014, in f\/inal form August 15, 2014; Published online August 19, 2014}

\Abstract{Maximal green sequences are particular sequences of mutations of quivers which were introduced by Keller in
the context of quantum dilogarithm identities and independently by Cecotti--C\'ordova--Vafa in the context of
supersymmetric gauge theory.
The existen\-ce of maximal green sequences for exceptional f\/inite mutation type quivers has been shown~by
Alim--Cecotti--C\'ordova--Espahbodi--Rastogi--Vafa except for the quiver $X_7$.
In this paper we show that the quiver $X_7$ does not have any maximal green sequences.
We also generalize the idea of the proof to give suf\/f\/icient conditions for the non-existence of maximal green sequences
for an arbitrary quiver.}


\Classification{15B36; 05C50}

\renewcommand{\thefootnote}{\arabic{footnote}}
\setcounter{footnote}{0}

\section{Introduction and main results}

Maximal green sequences are particular sequences of mutations of quivers.
They were used in~\cite{Ke2} to study the ref\/ined Donaldson--Thomas invariants and quantum dilogarithm identities.
Moreover, the same sequences appeared in theoretical physics where they yield the complete spectrum of a~BPS particle,
see~\cite[Section~4.2]{CCV}.
The existence of maximal green sequences for exceptional f\/inite mutation type quivers has been shown in~\cite{ACCERV}
except for the quiver~$X_7$.
In this paper, we show that the quiver~$X_7$ does not have any maximal green sequences.
We also give some general suf\/f\/icient conditions for the non-existence of maximal green sequences for an arbitrary
quiver.

To be more specif\/ic, we need some terminology.
Formally, a~quiver is a~pair $Q=(Q_0,Q_1)$ where $Q_0$ is a~f\/inite set of vertices and $Q_1$ is a~set of arrows between
them.
It is represented as a~directed graph with the set of vertices $Q_0$ and a~directed edge for each arrow.
We consider quivers with no loops or 2-cycles and represent a~quiver~$Q$ with vertices $1,\dots,n$,
by the uniquely associated skew-symmetric matrix $B=B^Q$ def\/ined as follows: the entry $B_{i,j}>0$ if and only if there are~$B_{i,j}$
many arrows from~$j$ to~$i$; if~$i$ and~$j$ are not connected to each other by an edge then $B_{i,j}=0$.
We will also consider more general skew-symmetrizable matrices: recall that an $n\times n$ integer matrix~$B$ is
skew-symmetrizable if there is a~diagonal matrix~$D$ with positive diagonal entries such that $DB$ is skew-symmetric.
To def\/ine the notion of a~green sequence, we consider pairs $(\mathbf{c}, B)$, where~$B$ is a~skew-symmetrizable integer
matrix and $\mathbf{c}=(\mathbf{c}_1,\dots,\mathbf{c}_n)$
such that each $\mathbf{c}_i=(c_1,\dots,c_n) \in \mathbb{Z}^n$ is non-zero.
Motivated by the structural theory of cluster algebras, we call such a~pair $(\mathbf{c}, B)$ a~$Y$-seed.
Then, for $k = 1, \dots, n$ and any~$Y$-seed $(\mathbf{c}, B)$ such that all entries of $ \mathbf{c}_k $ are
non-negative or all are non-positive, the \emph{$Y$-seed mutation} $\mu_k$ transforms $(\mathbf{c}, B)$ into
the~$Y$-seed $\mu_k(\mathbf{c}, B)=(\mathbf{c}', B')$ def\/ined as follows~\cite[equation~(5.9)]{CAIV}, where we use the
notation $[b]_+ = \max(b,0)$:
\begin{itemize}\itemsep=0pt
\item the entries of the exchange matrix $B'=(B'_{ij})$ are given~by
\begin{gather}
\label{eq:matrix-mutation}
B'_{ij} =
\begin{cases}
-B_{ij} & \text{if $i=k$ or $j=k$,}
\\
B_{ij} + [B_{ik}]_+ [B_{kj}]_+ - [-B_{ik}]_+ [-B_{kj}]_+ & \text{otherwise;}
\end{cases}
\end{gather}
\item the tuple $\mathbf{c}'=(\mathbf{c}_1',\dots,\mathbf{c}_n')$ is given~by
\begin{gather}
\label{eq:y-mutation}
\mathbf{c}'_i =
\begin{cases}
-\mathbf{c}_{i} & \text{if} \  i = k,
\\
\mathbf{c}_i+[\sgn(\mathbf{c}_k)B_{k,i}]_+\mathbf{c}_k & \text{if} \  i \neq k.
\end{cases}
\end{gather}
\end{itemize}
The matrix $B'$ is skew-symmetrizable with the same choice of~$D$.
We also use the notation $B' = \mu_k(B)$ (in~\eqref{eq:matrix-mutation}) and call the transformation $B \mapsto B'$ the
\emph{matrix mutation}.
This operation is involutive, so it def\/ines a~\emph{mutation-equivalence} relation on skew-symmetrizable matrices.

We use the~$Y$-seeds in association with the vertices of a~regular tree.
To be more precise, let~$\mathbb{T}_n$ be an \emph{$n$-regular tree} whose edges are labeled by the numbers $1, \dots,
n$, so that the~$n$ edges emanating from each vertex receive dif\/ferent labels.
We write $t \stackrel{k}{-} t'$ to indicate that vertices $t,t'\in\mathbb{T}_n$ are joined by an edge labeled by~$k$.
Let us f\/ix an initial seed at a~vertex $t_0$ in $\mathbb{T}_n$ and assign the (initial)~$Y$-seed $(\mathbf{c}_0,B_0)$,
where $\mathbf{c}_0$ is the tuple of standard basis.
This def\/ines a~\emph{$Y$-seed pattern} on $\mathbb{T}_n$, i.e.~an assignment of a~$Y$-seed $(\mathbf{c}_t, B_t)$ to
every vertex $t \in \mathbb{T}_n$, such that the seeds assigned to the endpoints of any edge $t \stackrel{k}{-} t'$ are
obtained from each other by the seed mutation~$\mu_k$; we call $(\mathbf{c}_t, B_t)$ a~\emph{$Y$-seed with respect to
the initial~$Y$-seed $(\mathbf{c}_0,B_0)$}.
We write:
\begin{gather*}
\mathbf{c}_t =\mathbf{c}= (\mathbf{c}_{1},\dots,\mathbf{c}_{n}),
\qquad
B_t=B = (B_{ij}).
\end{gather*}
We refer to~$B$ as the \emph{exchange matrix} and $\mathbf{c}$ as the \emph{$\mathbf{c}$-vector} tuple of the~$Y$-seed.
These vectors have the following \emph{sign coherence property}~\cite{DWZ2}:
\begin{gather}
\label{eq:C-sign-coherence}
\text{each vector $\mathbf{c}_{j}$ has either all entries nonnegative or all entries nonpositive.}
\end{gather}
Note that this property is conjectural if $B$ is a~general non-skew-symmetric (but skew-sym\-met\-ri\-zable) matrix.
It implies, in particular, that the {$Y$-seed mutation} in~\eqref{eq:y-mutation} is def\/ined for any $Y$-seed
$(\mathbf{c}_t, B_t)$, furthermore $\mathbf{c}_{t}$ is a~basis of $ \mathbb{Z}^n$~\cite[Proposition~1.3]{NZ}.
We also write $\mathbf{c}_{j}>0$ (resp.\
$\mathbf{c}_{j}<0$) if all entries are non-negative (resp.\
non-positive).

Now we can recall the notion of a~green sequence~\cite{BDP}:
\begin{definition}
Let $B_0$ be a~skew-symmetrizable $n\times n$ matrix.
A~\emph{green sequence for $B_0$} is a~sequence $\mathbf{i} = (i_1, \ldots, i_l)$ such that,   for any $1 \leq k \leq l$
with $(\mathbf{c},B)=\mu_{i_{k-1}} \circ \dots \circ \mu_{i_1}(\mathbf{c}_0,B_0)$, we have $\mathbf{c}_{i_{k}}>0$,
i.e.~each coordinate of $\mathbf{c}_{i_{k}}$ is greater than or equal to~$0$; here if $k=1$, then we take
$(\mathbf{c},B)=(\mathbf{c}_0,B_0)$.
\emph{A green sequence for a~quiver} is a~green sequence for the associated skew-symmetric matrix.

A green sequence $\mathbf{i} = (i_1, \ldots, i_l)$ is maximal if, for $(\mathbf{c},B)=\mu_{i_{l}} \circ \dots \circ
\mu_{i_1}(\mathbf{c}_0,B_0)$, we have $\mathbf{c}_k<0$ for all $k=1,\dots,n$.
\end{definition}

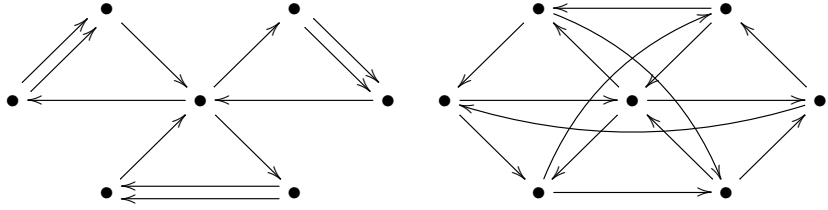
\begin{figure}[ht]
\centering
$
\begin{array}{cc}
\xymatrix{& \bullet\ar[rd] & & \bullet\ar@<.5ex>[rd]\ar@<-.5ex>[rd]
\\
\bullet\ar@<.5ex>[ru]\ar@<-.5ex>[ru] & & \bullet\ar[ll]\ar[ru]\ar[rd] & & \bullet\ar[ll]
\\
& \bullet\ar[ru] && \bullet\ar@<.5ex>[ll]\ar@<-.5ex>[ll] &} & \xymatrix{& \bullet\ar@/^1pc/[rrdd]\ar[ld] &&
\bullet\ar[ll]\ar[ld] &
\\
\bullet\ar[rd]\ar[rr] & & \bullet \ar[lu]\ar[rr]\ar[ld] & & \bullet \ar[lu]\ar@/^1pc/[llll]
\\
& \bullet\ar@/^1pc/[rruu]\ar[rr] & & \bullet\ar[lu]\ar[ru]}
\end{array}
$

\caption{Quivers which are mutation-equivalent to $X_7$; the f\/irst one is the quiver $X_7$, see~\cite{DO}.}
\label{fig:X7}
\end{figure}

In this paper, we study the maximal green sequences for the quivers which are mutation-equivalent to the quiver $X_7$
(Fig.~\ref{fig:X7}).
Our result is the following:

\begin{theorem}
\label{th:cyclic}
Suppose that~$Q$ is mutation-equivalent to the quiver $X_7$ $($so $Q$ is one of the quivers in Fig.~{\rm \ref{fig:X7})}.
Then~$Q$ does not have any maximal green sequences.
\end{theorem}

We prove the theorem using the following general statement, which can be easily checked to give a~suf\/f\/icient condition
for the non-existence of maximal green sequences:
\begin{proposition}
\label{prop:positive}
Let $B_0$ be a~skew-symmetrizable initial exchange matrix.
Suppose that there is a~vector $u>0$ such that, for any~$Y$-seed $(\mathbf{c},B)$ with respect to the initial seed
$(\mathbf{c}_0,B_0)$, the coordinates of~$u$ with respect to $\mathbf{c}$ are non-negative.
Then, under assumption~\eqref{eq:C-sign-coherence}, the matrix $B_0$ does not have any maximal green sequences.
\end{proposition}

We establish such a~vector for the quiver $X_7$:

\begin{proposition}
\label{prop:radical}
Suppose that $Q_0$ is a~quiver which is mutation-equivalent to $X_7$, so $Q_0$ is one of the quivers in
Fig.~{\rm \ref{fig:X7}}, and $B_0$ is the corresponding skew-symmetric matrix.
Let $u=(a_1,a_2,\dots,a_7)$ be the vector defined as follows:

\begin{quotation}
$(*)$~if $Q_0$ is the quiver $X_7$ $($so $Q$ is the f\/irst quiver in Fig.~{\rm \ref{fig:X7})}, then the coordinate
corresponding to the ``center'' is equal $2$, and the rest is equal to $1$; if $Q_0$ is not the quiver $X_7$ $($so $Q$ is the second quiver in Fig.~{\rm \ref{fig:X7})}, then all coordinates are equal to~$1$.
\end{quotation}

Then, for any~$Y$-seed $(\mathbf{c},B)$ with respect to the initial seed $(\mathbf{c}_0,B_0)$, the coordinates of~$u$
with respect to $\mathbf{c}$ is of the same form as in $(*)$.
In particular, the coordinates of~$u$ with respect to $\mathbf{c}$ are positive.

$($The vector $u$ is a~radical vector for $B_0$, i.e.~$B_0u=0$.
In fact, any radical vector for $B_0$ is a~multiple of $u.)$
\end{proposition}

We generalize this statement to an arbitrary quiver as follows:

\begin{theorem}
\label{th:radical}
Let $B_0$ be a~skew-symmetrizable initial exchange matrix and suppose that $u_0>0$ is a~radical vector for $B_0$,
i.e.~$ B_0u_0=0 $.
Suppose also that, for any~$Y$-seed $(\mathbf{c},B)$ with respect to the initial seed $(\mathbf{c}_0,B_0)$, the
coordinates of $u_0$ with respect to $\mathbf{c}$ are non-negative.
Then, under assumption~\eqref{eq:C-sign-coherence}, for any $B$ which is mutation-equivalent to $B_0$, the matrix $B$ does not have any maximal green sequences.
\end{theorem}

We prove our results in Section~\ref{sec:proof}. For related applications of maximal green sequences, we refer the reader to~\cite{CLS} and~\cite{S6}.

\section{Proofs of main results}
\label{sec:proof}

Let us f\/irst note how the coordinates of a~vector change under the mutation operation, which can be easily checked using
the def\/initions (assuming~\eqref{eq:C-sign-coherence}):

\begin{proposition}
\label{prop:coord_change}
Suppose that $(\mathbf{c},B)$ is a~$Y$-seed with respect to an initial~$Y$-seed.
Suppose also that the coordinate vector of~$u$ with respect to $\mathbf{c}$ is $(a_1,\dots,a_n)$.
Let $(\mathbf{c}', B')=\mu_k(\mathbf{c}, B)$ and $(a'_1,\dots,a'_n)$ be the coordinates of~$u$ with respect to
$\mathbf{c}'$.
Then $a_i=a'_i$ if $i \ne k$ and $a'_k=-a_k+\sum a_i[\sgn(\mathbf{c}_k)B_{k,i}]_+$, where the sum is over all $i \ne k$.
\end{proposition}

As we mentioned, in view of Proposition~\ref{prop:positive}, Theorem~\ref{th:cyclic} follows from
Proposition~\ref{prop:radical}.
To prove Proposition~\ref{prop:radical}, it is enough to show that the coordinates of the vector $u$ change as stated,
i.e.~show that if the coordinates of $u$ with respect to $\mathbf{c}$ are as in $(*)$, then for the~$Y$-seed
$(\mathbf{c}',B')=\mu_k(\mathbf{c},B)$, the coordinates with respect to $\mathbf{c}'$ are also of the form in $(*)$.
This can be checked easily using the formula in Proposition~\ref{prop:coord_change}.

To prove Theorem~\ref{th:radical}, let us f\/irst note the following property of the coordinates of the radical vectors:

\begin{lemma}
\label{lem:rad_coord_change}
Suppose that $(\mathbf{c},B)$ is a~$Y$-seed with respect to an initial~$Y$-seed $(\mathbf{c}_0,B_0)$ and $u_0$ is
a~radical vector for $B_0$.
Suppose that the coordinate vector of $u_0$ with respect to $\mathbf{c}$ is $(a_1,\dots,a_n)$.
Then, for any index $k$, we have the following:
\begin{gather*}
\sum a_i[\sgn(\mathbf{c}_k)B_{k,i}]_+=\sum a_i[-\sgn(\mathbf{c}_k)B_{k,i}]_+,
\end{gather*}
where the sum is over all $i \ne k$.

In particular, for radical vectors, the formula in Proposition~{\rm \ref{prop:coord_change}} that describe the change of
coordinates under mutation depends only on the exchange matrix, not on the $\mathbf{c}$-vectors.
\end{lemma}

To prove the lemma, suppose that $ D=\operatorname{diag}(d_1,\dots,d_n) $ is a~skew-symmetrizing matrix for~$B_0$, so it is also
skew-symmetrizing for~$B$, so $ DB=C $ is skew-symmetric, i.e.~$ C_{i,k}=d_iB_{i,k}=-d_kB_{k,i}=-C_{k,i} $ for all~$i$,~$k$.
Let $ u=(a_1,\dots,a_n) $.
Then $u$ is a~radical vector for $B$, so it is also a~radical vector for $C=DB$, i.e.~$ Cu=0 $, which means that
for any index $k$, we have $ \sum a_i[\sgn(\mathbf{c}_k)C_{k,i}]_+=\sum a_i[\sgn(\mathbf{c}_k)C_{i,k}]_+$, which is
equal to $\sum a_i[-\sgn(\mathbf{c}_k)C_{k,i}]_+$, where the sum is over all $i \ne k$.
Then, writing $ C_{k,i}= d_kB_{k,i}$, we have
\begin{gather*}
\sum a_i[\sgn(\mathbf{c}_k)d_kB_{k,i}]_+=\sum a_i[-\sgn(\mathbf{c}_k)d_kB_{k,i}]_+.
\end{gather*}
Dividing both sides by $d_k>0$, we obtain the lemma.

 We will also need the following property of the radical vectors:

\begin{lemma}
\label{lem:rad_coord_change2}
In the set-up of Theorem~{\rm \ref{th:radical}}, let $u$ denote the vector which represents $u_0$ with respect to the
basis $\mathbf{c}$.
Then $u$ is a~radical vector for $B$, i.e.~$Bu=0$.
\end{lemma}

To prove the lemma, let us note that $u$ can be obtained from $u_0$ by applying the formula in
Proposition~\ref{prop:coord_change} along with the mutations.
Thus, to prove the lemma, it is enough to show that, for any $k=1,\dots,n$, we have the following:
\begin{quotation}
$(**)$~if $u=(a_1,\dots,a_n)$ is a~radical vector for $B$, then $u'$ is a~radical vector for $B'=\mu_k(B)$, i.e.~$B'u'=0$,
where $u'=(a'_1,\dots,a'_n)$ is the vector as in Proposition~\ref{prop:coord_change}.
\end{quotation}

To show $(**)$, we write $B'$ in matrix notation as follows~\cite[Lemma~3.2]{CAIII}: for $\epsilon=\sgn(\mathbf{c}_k)$,  we have
\begin{gather*}
B' = (J_{n,k} + E_k)  B  (J_{n,k} + F_k),
\end{gather*}
where
\begin{itemize}\itemsep=0pt
\item $J_{n,k}$ denotes the diagonal $n\times n$ matrix whose diagonal entries are all~$1$'s, except for $-1$ in
the~$k$th position;
\item $E_k$   is the $n\times n$ matrix whose only nonzero entries are   $e_{ik} = [-\varepsilon b_{ik}]_+$;
\item $F_k$   is the $n\times n$ matrix whose only nonzero entries are $f_{kj} = [\varepsilon b_{kj}]_+$.
\end{itemize}
It follows from a~direct check that $(J_{n,k} + F_k)u'=u$.
Then $ B'u'= (J_{n,k} + E_k) B (J_{n,k} + F_k)u'= (J_{n,k} + E_k) Bu=(J_{n,k} + E_k) 0=0$.
This completes the proof of the lemma.

Let us now prove Theorem~\ref{th:radical}.
For this, we f\/irst consider the $Y$-seed pattern def\/ined by the initial $Y$-seed $(\mathbf{c}_0,B_0)$ at the initial
vertex $t_0$.
Let us suppose that $t_1$ is a~vertex such that the corresponding $Y$ seed $(\mathbf{c},B)$ has the exchange matrix~$B$.
Then we can consider the $Y$-seed pattern def\/ined by the initial $Y$-seed $(\mathbf{c}_0,B)$ at the initial vertex
$t_1$ (where $\mathbf{c}_0$ is the tuple of standard basis).
Then we have the following: for any f\/ixed vertex $t$ of the $n$-regular tree $\mathbb{T}_n$, the exchange matrices
of the $Y$-seeds assigned by these patterns coincide because the pattern is formed by mutating at the labels
of the $n$-regular tree $\mathbb{T}_n$ and mutation is an involutive operation on matrices; let us denote these seeds~by
$(\mathbf{c}',B')$ and $(\mathbf{c}'',B')$ respectively.

On the other hand, let $u$ denote the vector which represents~$u_0$ with respect to the basis~$\mathbf{c}$, which
can be obtained by applying the formula in Proposition~\ref{prop:coord_change} along with the mutations.
Then~$u$ is a~radical vector for~$B$, i.e.~$Bu=0$ (Lemma~\ref{lem:rad_coord_change2}).
Furthermore, the coordinates of the vectors $u_0$ and $u$ with respect to the bases~$\mathbf{c}'$ and~$\mathbf{c}''$
respectively will coincide by Lemma~\ref{lem:rad_coord_change} (which says that for radical vectors the formula in
Proposition~\ref{prop:coord_change} depends only on the exchange matrices, not on the $ \mathbf{c}$-vectors).
In particular, the coordinates of~$u$ with respect to any basis of~$ \mathbf{c} $-vectors are non-negative.
Thus, by Proposition~\ref{prop:positive}, the matrix $B$ does not have any maximal green sequences.
This completes the proof.

\subsection*{Acknowledgements}

The author's research was supported in part by the Scientif\/ic and Technological Research Council of Turkey (TUBITAK)
grant \#~113F138.
The author also thanks Christof\/f Geiss for drawing his attention to the paper~\cite{ACCERV} by presenting its results at
the Workshop on Hall and Cluster Algebras in CRM, University of Montreal.
He also thanks the organizers for organizing the conference.

\pdfbookmark[1]{References}{ref}
\LastPageEnding

\end{document}